\documentclass[12pt]{amsart}
\usepackage{amsmath}
\usepackage{amssymb}
\usepackage{mathrsfs}
\newtheorem{theorem}{Theorem}
\newtheorem{lem}[theorem]{Lemma}
\newtheorem{prop}[theorem]{Proposition}
\newtheorem{theoremb}{Theorem}
\newtheorem{dfn}[theoremb]{Definition\!\!}
\newtheorem{rk}[theoremb]{Remark\!\!}

\newcommand\bib[1]{\bibitem[#1]{#1}}

\newcommand\p{\partial}
\newcommand\D{{\mathcal D}}
\newcommand\z{\sigma}
\renewcommand\k{\varkappa}
\newcommand\R{{\mathbb R}}
\newcommand\op[1]{\mathop{\rm #1}\nolimits}
\newcommand\Cc{{\let\mathcal\mathscr\mathcal C}}
\DeclareFontFamily{U}{wncyr}{}
\DeclareFontShape{U}{wncyr}{m}{it}{%
  <5><6><7><8><9>gen*wncyi%
  <10><10.95><12><14.4><17.28><20.74><24.88>wncyi10}{}
\DeclareSymbolFont{MathRussLetters}{U}{wncyr}{m}{it}
\DeclareMathSymbol{\re}{\mathalpha}{MathRussLetters}{3}

\begin{document}

 \title{Anomaly of linearization \\ and auxiliary integrals.}
 \author{Boris Kruglikov}
 \address{Institute of Mathematics and Statistics, University of Troms\o, Troms\o\ 90-37, Norway.}
 \email{kruglikov@math.uit.no}
 \date{}
 \maketitle

 \vspace{-14.5pt}
 \begin{abstract}
In  this note we discuss some formal properties of universal linearization operator,
relate this to brackets of non-linear differential operators and discuss application
to the calculus of auxiliary integrals, used in compatibility reductions of PDEs.
 \end{abstract}

\section*{Introduction}\label{S0}

Commutator $[\Delta,\nabla]$ of linear differential operators $\Delta,\nabla\in\op{Diff}(\pi,\pi)$ in the context of non-linear operators $F,G\in\op{diff}(\pi,\pi)$ is up-graded to the higher Jacobi bracket $\{F,G\}$, which plays the same role in compatibility investigations and symmetry calculus.\footnote{MSC numbers: 35A27, 58A20; 58J70, 35A30.\\
Keywords: Linearization, evolutionary differentiation, compatibility, differential constraint, symmetry, reduction, Jacobi bracket, multi-bracket.}

The linearization operator relates non-linear operators on a bundle $\pi$ with linear operators on the same bundle, whose coefficients should be however smooth functions on the space of infinite jets. The latter space is the algebra of $\Cc$-differential operators and we get the map
 $$
\ell:\op{diff}(\pi,\pi)\to\Cc\op{Diff}(\pi,\pi)=C^\infty(J^\infty\pi)\otimes_{C^\infty(M)}\mathop{\rm Diff}(\pi,\pi),
 $$
defined by the formula \cite{KLV}
 $$
\ell_{F}(s)h=\tfrac{d}{dt}F(s+th)|_{t=0},\qquad\ F\in\op{diff}(\pi,\pi),\quad s,h\in C^\infty(\pi).
 $$
However it does not respect the commutator:
 $$
[\ell_F,\ell_G]\ne\ell_{\{F,G\}}.
 $$

\hspace{-13.5pt}{\bf Example:} Consider the scalar differential operators on $\R$, so that $\pi=\mathbf{1}$
and $J^\infty(\pi)=\R^\infty(x,u,p=p_1,p_2,\dots)$. Choose
 $$
F=p^2,G=p+c\cdot x;\quad \{F,G\}=2c\,p\,\Longrightarrow\,\ell_{\{F,G\}}=2c\,\mathcal{D}_x.
 $$
If we commute $\ell_F=2p\,\mathcal{D}_x$ and $\ell_G=\mathcal{D}_x$, we get: $[\ell_F,\ell_G]=-2p_2\,\mathcal{D}_x$, so that we observe an anomaly.

There are two reasons for this. The first is that the operator of linearization disregards non-homogeneous linear terms, which are important for the Jacobi bracket. The second is the non-linearity itself.

The goal of this note is to discuss reasons and consequences of this anomaly (this also plays a significant role in investigation of coverings and non-local calculus \cite{KKV}).

\textsc{Acknowledgement.}
The results were obtained and systematized during the research stay in Max Planck Institute for Mathematics in the Sciences, Leipzig, in April-May 2007.

\section{Anomaly via Hessian}\label{S1}

The Jacobi bracket of non-linear operators $F,G\in\op{diff}(\pi,\pi)$ is expressed via linearization as follows:
 $$
\{F,G\}=\ell_FG-\ell_GF.
 $$

We also consider the evolutionary operators defined by duality:
 $$
\re_FG=\ell_GF.
 $$
Since $\ell_G$ is a derivation in $G$, $\re_F$ is a derivation (satisfies the Leibniz rule) and their union can be treated as the module of vector fields.
These operators have no anomaly, i.e. the map $\re:C^\infty(J^\infty\pi)\to\op{Vect}(J^\infty\pi)$ is an anti-homomorphism:
 $$
[\re_F,\re_G]=-\re_{\{F,G\}}.
 $$
This instantly implies Jacobi identity for the bracket $\{F,G\}$, so that $\bigl(\op{diff}(\pi,\pi),\{,\}\bigr)$ is a Lie algebra \cite{KLV}.

The operators of universal linearization and evolutionary differentiation do not commute and this leads to the following

 \begin{dfn}
The Hessian operator
$\op{diff}(\pi,\pi)\times\op{diff}(\pi,\pi)\to\Cc\op{Diff}(\pi,\pi)$ is defined by the formula
 $$
\op{Hess}_FG=[\re_G,\ell_F].
 $$
 \end{dfn}
We will also write $\op{Hess}_F(G,H)=\op{Hess}_FG(H)$ for $F,G,H\in\op{diff}(\pi,\pi)$ and note that $\op{Hess}_F\equiv0$ for linear operators $F$, because in this case $\ell_F=F$, which reduces the claim to the commutation of left and right multiplications.

Next we note that the Hessian $\op{Hess}_F$ is symmetric:

 \begin{lem}
$\op{Hess}_F(G,H)=\op{Hess}_F(H,G)$.
 \end{lem}

Indeed:
 $$
\op{Hess}_F(G,H)=\re_G\ell_FH-\ell_F\re_GH=\re_G\re_HF-\ell_F\ell_HG,
 $$
so that
 \begin{multline*}
\op{Hess}_F(G,H)-\op{Hess}_F(H,G)=[\re_G,\re_H]F-\ell_F\{H,G\}\\
=-\re_{\{G,H\}}F-\ell_F\{H,G\}=0.
 \end{multline*}

Now we can express the anomaly of linearization via the Hessian:

 \begin{prop}\label{pro2}
$[\ell_F,\ell_G]-\ell_{\{F,G\}}=\op{Hess}_GF-\op{Hess}_FG$.
 \end{prop}

Indeed we have:

 \begin{multline*}
[\ell_F,\ell_G]H=\ell_F\re_HG-\ell_G\re_HF\\
=\re_H(\ell_FG-\ell_GF)-\op{Hess}_F(H,G)+\op{Hess}_G(H,F)\\
=\re_H\{F,G\}-\op{Hess}_F(G,H)+\op{Hess}_G(F,H)\\
=\ell_{\{F,G\}}H+(\op{Hess}_GF-\op{Hess}_FG)H.
 \end{multline*}

Finally let us express the Leibniz identity for non-linear operators and the Jacobi bracket.
For linear operators it is well-known, but for non-linear ones there's an anomaly:

 \begin{prop}\label{pro3}
$\{F,\ell_GH\}=\ell_{\{F,G\}}H+\ell_G\{F,H\}-\op{Hess}_F(G,H)$.
 \end{prop}

This is obtained as follows:

 \begin{multline*}
\{F,\ell_GH\}=\ell_F\ell_GH-\re_F\ell_GH\\
=[\ell_F,\ell_G]H+\ell_G(\ell_F-\re_F)H-\op{Hess}_G(F,H)\\
=\ell_{\{F,G\}}H+\ell_G\{F,H\}-\op{Hess}_F(G,H).
 \end{multline*}

\section{Coordinate expressions}\label{S2}

A local coordinate system $(x^i,u^j)$ on $\pi$ induces the canonical coordinates $(x^i,p^j_\z)$ on the space $J^\infty\pi$, where $\z=(i_1,\dots,i_n)$ is a multi-index of length $|\z|=i_1+\dots+i_n$. The operator of total derivative of multi-order $\z$ (and order $|\z|$) is $\D_\z=\D_1^{i_1}\cdots\D_n^{i_n}$, where $\D_i=\p_{x^i}+\sum p_{\tau+1_i}^j\p_{p_\tau^j}$.

The linearization of $F=(F_1,\dots,F_r)$ is $\ell_F=(\ell(F_1),\dots,\ell(F_r))$
with
 $$
\ell(F_i)=\sum (\p_{p_\z^j}F_i)\cdot \D_\z^{[j]},
 $$
where $\D_\z^{[j]}$ denotes the operator $\D_\z$ applied to the
$j$-th component of the section from $C^\infty(\pi)$.

The $i$-th component of the evolutionary differentiation $\re_G$ corresponding to $G=(G_1,\dots,G_n)$ equals
 $$
\re_G^i=\sum(\D_\z G_j)\cdot\p_{p_\z^j}{}^{[i]},
 $$
where $\p_{p_\z^j}{}^{[i]}$ denotes the operator $\p_{p_\z^j}$
applied to the $i$-th component of the section from $C^\infty(\pi)$.

Then $i$-th components of the Jacobi bracket is given by
 $$
\{F,G\}_i=\sum \bigl(\D_\z(G_j)\cdot\p_{p_\z^j}F_i-
\D_\z(F_j)\cdot\p_{p_\z^j}G_i\bigr).
 $$

These formulas are known \cite{KLV}. It is instructive to demonstrate the Jacobi identity in coordinates.
For this we need the following assertion.
 \begin{lem}
In canonical coordinates on $J^\infty\pi$:
 $$
\p_{p_\z^i}\D_\tau=\sum\D_{\tau-\k}\p_{p_{\z-\k}^i}
 $$
(the difference of multi-indices $\z-\k$ is defined whenever $\k\subset\z$), the summation is by $\k$ counted with multiplicity. More generally for vector differential operators if $\D^{[j]}_\z$ is the operator $\D_\z$ acting on the $j$-th component, then the above formula holds true for such specification.
 \end{lem}

This follows from iteration of the formula $[\p_{p^j_\z},\D_i]=\p_{p^j_{\z-1_i}}$. Thus
 \begin{multline*}
 \{F,\{G,H\}\}=\sum F_{p_\z}\D_{\z-\k}(G_{p_\tau})\,\D_{\tau+\k}(H)- F_{p_\z}\D_{\z-\k}(H_{p_\tau})\,\D_{\tau+\k}(G)\\
-G_{p_\z p_\tau}\D_\tau(H)\D_\z(F)+H_{p_\z p_\tau}\D_\tau(G)\D_\z(F)\\
-(G_{p_\z}\D_{\z-\k}(H_{p_{\tau-\k}})-H_{p_\z}\D_{\z-\k}(G_{p_{\tau-\k}}))\D_\tau(F),
 \end{multline*}
which yields $\sum_\text{cyclic}\{F,\{G,H\}\}=0$.

Now we write the Hessian:
 $$
\op{Hess}_F(G,H)=\sum F_{p_\z p_\tau}\D_{\z}G\cdot\D_\tau H,
 $$
and its symmetry in $G,H$ and vanishing for linear $F$ is obvious.

The compensated Leibniz formula can be written as follows:

 \begin{multline*}
\{F,\ell_GH\}-\ell_{\{F,G\}}H-\ell_G\{F,H\}=\\
\sum F_{p_\z}\D_{\z-\k}(G_{p_\tau})\,\D_{\tau+\k}(H)-
(G_{p_\z p_\tau}\D_\tau(H)\D_\z(F)+G_{p_\tau}\p_{p_\z}\D_\tau(H))\D_\z(F)\\
-(F_{p_\z p_\tau}\D_\z(G)+F_{p_\z}\p_{p_\tau}\D_\z(G))\D_\tau(H)+
(G_{p_\z p_\tau}\D_\z(F)+G_{p_\z}\p_{p_\tau}\D_\z(F))\D_\tau(H)\\
-G_{p_\z}\,(\D_{\z-\k}(F_{p_\tau})\D_{\tau+\k}(H)-\D_{\z-\k}(H_{p_\tau})\D_{\tau+\k}(F))
 =-\op{Hess}_F(G,H)
 \end{multline*}
and the anomaly in commuting linearizations is:
 \begin{multline*}
[\ell_F,\ell_G]-\ell_{\{F,G\}}=\\
\sum F_{p_\z}\D_{\z-\k}(G_{p_\tau})\,\D_{\tau+\k}(H)- G_{p_\z}\D_{\z-\k}(F_{p_\tau})\,\D_{\tau+\k}(H)\\
-(F_{p_\z p_\tau}\D_\z(G)+F_{p_\z}\p_{p_\tau}\D_\z(G))\D_\tau(H)+
(G_{p_\z p_\tau}\D_\z(F)+G_{p_\z}\p_{p_\tau}\D_\z(F))\D_\tau(H)\\
 =\op{Hess}_G(F,H)-\op{Hess}_F(G,H).
 \end{multline*}

This gives an alternative proof of Propositions \ref{pro3} and \ref{pro2}.

\section{Auxiliary integrals}\label{S3}

 \begin{dfn}
An operator $G\in\op{diff}(\pi,\pi)$ is called an auxiliary integral for $F\in\op{diff}(\pi,\pi)$ if
 $$
\{F,G\}=\ell_\lambda F+\ell_\mu G
 $$
for some operators $\lambda\in\op{diff}(\pi,\pi)$ and $\mu\not\in\Cc\op{Diff}(\pi,\pi)\cdot F\setminus\{0\}$. The set of such $G$ is denoted by $\op{Aux}(F)$.
 \end{dfn}

It is better to denote $\op{Aux}_\mu(F)$ the space of $G$ satisfying the above formula with some fixed $\mu\in\op{diff}(\pi,\pi)$, because it is a vector space. Then $\op{Aux}(F)=\cup_\mu\op{Aux}_\mu(F)$. We can assume $\op{ord}(\mu)<\op{ord}(F)$ for scalar operators, i.e. $\op{rank}\pi=1$.

With certain non-degeneracy condition for the symbols of $F,G$ the following statement holds:

 \begin{theorem}\label{thm5}
A non-linear differential operator $G$ is an auxiliary integral for another operator $F$ iff the system $F=0,G=0$ is compatible (formally integrable).
 \end{theorem}

The generic position condition for the symbols of $F,G$ is essential. If $\pi={\bf1}$ is the trivial one-dimensional bundle, this condition is just the transversality of the characteristic varieties $\op{Char}^\mathbb{C}(F)$ and $\op{Char}^\mathbb{C}(G)$ in the bundle $\mathbb{P}^\mathbb{C}T^*M$ (after pull-back to the joint system $F=G=0$ in jets); in this form it is a particular form of the statement proved in \cite{KL$_2$}. For $\op{rank}\pi>1$ the condition is more delicate and will be presented elsewhere.

Notice that $\op{Aux}_0(F)=\op{Sym}(F)$ is the space of symmetries of $F$. This is a Lie algebra with respect to the Jacobi bracket. It can be represented as a union of spaces
 $$
\op{Sym}_\theta(F)=\{H\,:\,\ell_FH=\ell_{\theta+H}F\},\qquad \theta\in\op{diff}(\pi,\pi),
 $$
which are modules over $\op{Sym}_0(F)$. More generally we have the graded group: $\op{Sym}_{\theta'}(F)+\op{Sym}_{\theta''}(F)\subset\op{Sym}_{\theta'+\theta''}(F)$

Let us assume $G\in\op{Aux}_\mu(F)$, $H\in\op{Sym}_\theta(F)$, i.e.
 $$
\{F,G\}=\ell_\lambda F+\ell_\mu G,\qquad \{F,H\}=\ell_\theta F.
 $$
Then denoting $\op{ad}_H=\{H,\cdot\}=\ell_H-\re_H$ we get:

 \begin{multline*}
 \op{ad}_F\{G,H\}= \{\op{ad}_FG,H\}+\{G,\op{ad}_FH\}\\
= -\{H,\ell_\lambda F+\ell_\mu G\}+\{G,\ell_\theta F\}\\
=\ell_{\{\lambda,H\}}F+\ell_{\lambda}\{F,H\}+\op{Hess}_H(\lambda,F)
+\ell_{\{\mu,H\}}G+\ell_{\mu}\{G,H\}+\op{Hess}_H(\mu,G)\\
-\ell_{\{\theta,G\}}F-\ell_\theta\{F,G\}-\op{Hess}_G(\theta,F)\\
= (\ell_{\{\lambda,H\}}+[\ell_\lambda,\ell_\theta]-\ell_{\{\theta,G\}}+\op{Hess}_H\lambda-
\op{Hess}_G\theta)F+\ell_\mu\{G,H\}\\
 +(\ell_{\{\mu,H\}}-\ell_\theta\ell_\mu+\op{Hess}_H\mu)G.
 \end{multline*}

Thus $\{G,H\}$ is an auxiliary integral for $F$ if $\ell_\theta\ell_\mu=\ell_{\{\mu,H\}}+\op{Hess}_H\mu$ (the "iff" condition means the difference annihilates $G$), which can be written as
 $$
\mu\in\op{Ker}[(\ell_\theta+\ell_{\op{ad}_H}-\op{Hess}_H)\circ\ell\,].
 $$

Such a pair $\theta\in\op{sym}^*(F)=\op{Sym}(F)/\op{Sym}_0(F)$, $H\in\op{Sym}_\theta(F)$ determines the action of the second component
 $$
\op{ad}_H:\op{Aux}_\mu(F)\to\op{Aux}_\mu(F).
 $$

Also since
 \begin{multline*}
\ell_{\{\mu,H\}}G=\re_G\{\mu,H\}=\re_G\re_H(\mu)-\re_G\ell_H(\mu)
=(\re_H-\ell_H)\re_G(\mu)\\
-\re_{\{G,H\}}\mu-\op{Hess}_H(G,\mu)=
-\op{ad}_H\re_G(\mu)-\op{Hess}_H(\mu,G)-\ell_\mu\{G,H\},
 \end{multline*}
we have:
 $$
\ell_\mu\{G,H\}+(\ell_{\{\mu,H\}}-\ell_\theta\ell_\mu+\op{Hess}_H\mu)G=
-(\op{ad}_H+\ell_\theta)\ell_\mu G.
 $$
Thus if $H\in\op{Sym}_\theta(F)$, i.e. $(\op{ad}_H+\ell_\theta)F=0$, and $\mu\in\op{Ker}[(\op{ad}_H+\ell_\theta)\circ\ell\,]$, i.e.
$(\op{ad}_H+\ell_\theta)\ell_\mu=0$, then
 $$
\op{ad}_H:\op{Aux}_\mu(F)\to\op{Sym}(F).
 $$

\section{Symmetries and compatibility}\label{S4}

It has been a common belief that if $G\in\op{Sym}(F)$, then the system $F=0,G=0$ is compatible, which forms the base of investigation for auto-model solutions. This is however not always true.

\smallskip

\hspace{-13.5pt}{\bf Example:} Let $F,G$ be two linear diagonal operators with constant coefficients. Then $\{F,G\}=0$ (in this case the Jacobi bracket is the standard commutator), so that $G$ is a symmetry of $F$. However the system $F=0,G=0$ is usually incompatible: for generic $F,G$ of the considered type the only solution will be the trivial zero vector-function.

More complicated non-diagonal operators are possible, but it would be better to consider non-homogeneous linear operators. Then if the coefficients are constant and generic, the linear matrix part commute, but the system $F=0,G=0$ may have no solutions at all.

For instance if we take
 \begin{gather*}
F=\begin{bmatrix}(\D_x^2-\D_y) & 0 \\ 0 & (\D_x\D_y+1) \end{bmatrix}\cdot \begin{bmatrix}
u \\ v\end{bmatrix}-\begin{bmatrix}1 \\ 0\end{bmatrix},\\
G=\begin{bmatrix}(\D_x\D_y-1) & 0 \\ 0 & (\D_y^2-\D_x) \end{bmatrix}\cdot \begin{bmatrix}
u \\ v\end{bmatrix}+\begin{bmatrix}0 \\ 1\end{bmatrix},
 \end{gather*}
then $\{F,G\}=0$, so that $G\in\op{Sym}(F)$, while the system $F=0,G=0$ is not compatible, and moreover its solutions space is empty.

\smallskip

Thus the flow $u_t=G(u)$ on the equation $F=0$ has no fixed points (no auto-model solutions).
Here $t$ is an additional variable ($x$ is the base multi-variable for PDEs $F=0$ and $G=0$), so that $G\in\op{Sym}(F)$ can be expressed as compatibility of the system
 $$F(u)=0,\quad u_t=G(u),$$
while symmetric solutions correspond to the stationary case $u_t=0$, i.e. compatibility of the system $F(u)=0,G(u)=0$\footnote{I am grateful to S.Igonin and A.Verbovetsky for an enlightening discussion about the results of \cite{KL$_2$,KL$_3$} and the symmetry condition.}.

However if the non-degeneracy condition assumed in Theorem \ref{thm5} is satisfied, then auto-model (or invariant) solutions exist in abundance, namely they have the required functional dimension and rank as Hilbert polynomial (or Cartan test \cite{C}) predicts, see \cite{KL$_4$}.

 \begin{rk}
Symmetric solutions are the stationary points of the evolutionary fields and they are similar to the fixed points of smooth vector fields on $\R^n$, which must exist provided the vector field is Morse at infinity. The non-degeneracy condition plays a similar role.
 \end{rk}

Many examples of auto-model solutions and their generalizations can be found in \cite{BK,Ol,Ov}, non-local analogs use the same technique and similar theory \cite{KLV,KK,KKV}.

Compatible systems correspond to reductions of PDEs and are sometimes called conditional symmetries by analogy with finite-dimensional integrable systems on one isoenergetic surface \cite{FZ}. But the rigorous result must rely on certain general position property for the symbol of differential operators, otherwise it can turn wrong \cite{KL$_2$,KL$_3$}. The method based on this approach makes specification of the general idea of differential constraint and is described in \cite{KL$_1$}.

\section{Conclusion}\label{S5}

In this note we described the higher-jets calculus corresponding to symmetries and compatible constraints, basing on the Jacobi brackets. Another approach to integrability of vector systems is given by minimal overdetermination and it uses multi-brackets of differential operators
 $$
\{\cdots\}:\Lambda^{m+1}\op{diff}(m\cdot{\bf1},{\bf1})\to\op{diff}(m\cdot{\bf1},{\bf1})
 $$
introduced in \cite{KL$_3$}, which are governed by the non-commutative Pl\"ucker identity. 

Following this approach a minimal generalization of symmetry for $F=(F_1,\dots,F_m)\in\op{diff}(\pi,\pi)$ with $\pi=m\cdot{\bf1}$ is such $G\in\op{diff}(\pi,{\bf1})$ that
 $$
\{F_1,\dots,F_m,G\}=\ell_{\theta_1}F_1+\dots+\ell_{\theta_m}F_m.
 $$
With certain non-degeneracy assumption \cite{KL$_3$} this implies that the overdetermined system $F=0,G=0$ is compatible (formally integrable).

A more advanced algebraic technique would yield another higher-jets calculus producing anomaly that manifests in non-vanishing of the expression
 $$
\{\ell_{F_1},\cdots,\ell_{F_{m+1}}\}-\ell_{\{F_1,\cdots,F_{m+1}\}}.
 $$
Implications for vector auxiliary integrals and generalized Lagrange-Charpit method follow the same scheme.


\end{document}